\documentclass[12pt]{amsart}
\usepackage{graphicx,amsmath,amssymb,extsizes,fourier,xcolor}
\usepackage[font=small]{caption}
\usepackage{accents}
\usepackage{float}
\usepackage{soul}
\usepackage[shortlabels]{enumitem}

\newcommand{\doubletilde}[1]{\widetilde{\mathop{\widetilde{#1}}}}

\begin{document}

\title{An elementary derivation of 3-cycles for a quadratic map}

\subjclass[2020]{Primary 39A23, 39A30; Secondary 39A33, 37E05}

\keywords{quadratic map, logistic map, 3-cycle, stability, conjugacy of maps, chaos}

\author{\'Arp\'ad B\'enyi}
\address{Department of Mathematics, Western Washington University,
516 High St, Bellingham, WA 98225, USA}
\email{benyia@wwu.edu}

\author{Ioan Ca\c{s}u}
\address{Department of Mathematics, West University of Timisoara, Bd. Vasile P\^arvan, nr. 4, 300223 Timi\c{s}oara, Romania}
\email{ioan.casu@e-uvt.ro}
\thanks{Corresponding author: Ioan Ca\c{s}u \texttt{ioan.casu@e-uvt.ro}.}

\date{}

\begin{abstract}
We present an elementary derivation of the period-three cycles for the real quadratic map $x\mapsto x^2+c$, a fundamental model in one-dimensional discrete dynamics. Using symmetric polynomials, we obtain a complete algebraic characterization of 3-cycles and determine explicit conditions for their existence and stability, without reliance on computer algebra. Through conjugacy with the logistic map, we recover the classical threshold values of the logistic parameter corresponding to the emergence and loss of stability of the 3-cycle. Our methodology outlines a transparent and algebraically grounded route to understanding the onset of chaos in quadratic and logistic dynamics. 
\end{abstract}

\maketitle

\section{Introduction}

The work of Li and Yorke \cite{LY}, building on Sharkovski\u{\i}'s theorem \cite{Sh}, provides one of the simplest, and arguably most famous, criteria for the onset of chaos in one-dimensional dynamics:

\centerline{\emph{If $f(x)$ has a period-three cycle, then $x_{n+1}=f(x_n)$ is chaotic.}}

\smallskip

As May observed in his seminal paper \cite{Ma}, many simple mathematical models can exhibit remarkably complex dynamics. Among these, the logistic and quadratic maps stand as paradigmatic examples of chaos. In what follows, we focus on the real version of ``the'' complex quadratic map that defines the Mandelbrot set \cite{Man}; namely, for a given parameter $c\in\mathbb R$, we let
\begin{equation}
    \label{map}
    f_c: \mathbb{R}\to \mathbb{R},~~f_c(x)=x^2+c;
\end{equation}
see, for example, Hasselblatt-Katok's book \cite[Chapter 11]{HK} for a thorough study of quadratic maps. 

The modest aim of this work is to present an elementary approach to the study of the period-three case for the quadratic map, viewed as a natural gateway to chaos. More precisely, we show that a complete 3-cycle analysis can be done via symmetric polynomials and without computer algebra. Our hope is that the route to the period-three case exposed here allows the reader to better grasp the deep connection between algebraic manipulations and dynamical phenomena. Through the use of conjugation, we also recover several classical results known for the logistic map. It is worthwhile mentioning that, in light of the conjugacy relation \eqref{r-c}, it seems more natural to first resolve the quadratic map case and then deduce the logistic map case from it. Indeed, while each logistic map parameter uniquely determines its corresponding quadratic map parameter, the converse holds only when the quadratic map parameter lies within the range of the quadratic function in \eqref{r-c}.

\section{Existence and stability of 3-cycles for the quadratic map}
For simplicity, in what follows we prefer to write $f$ instead of $f_c$. The first goal of this section is to determine the values of the parameter $c$ such that there exists a 3-cycle $(x_1,x_2,x_3)$ for the map \eqref{map}. Namely, we investigate the existence of distinct real numbers $x_1, x_2, x_3$ such that
\begin{align}
  &  x_2=x_1^2+c=f(x_1), \label{eq1}\\
   & x_3=x_2^2+c=f(x_2), \label{eq2}\\
  &  x_1=x_3^2+c=f(x_3). \label{eq3}
\end{align}
The symmetric polynomials in the variables $x_1, x_2, x_3$ will be denoted by
$$s_1=x_1+x_2+x_3,~~~~s_2=x_1x_2+x_2x_3+x_3x_1, ~~~~\text{and}~~~~s_3=x_1x_2x_3.$$
Adding equations \eqref{eq1}-\eqref{eq3} we obtain
\begin{equation}
    \label{s-1}
    s_1=s_1^2-2s_2+3c.
\end{equation}
Multiplying equations \eqref{eq1}-\eqref{eq3} with $x_1,x_2$ and $x_3$ respectively, and then adding them, we obtain
\begin{equation}
    \label{s-2}
 s_2=3s_3+s_1(s_1^2-3s_2)+cs_1=3s_3+s_1^3-3s_1s_2+cs_1.   
\end{equation}
Also,
\begin{equation*}
\begin{aligned}
    s_2^2&=(x_2x_3+x_3x_1+x_1x_3)^2=\sum x_2^2x_3^2+2s_1s_3\\
    &=\sum (x_3-c)(x_1-c)+2s_1s_3=s_2-2cs_1+3c^2+2s_1s_3,
\end{aligned}
\end{equation*}
that is,
\begin{equation}
    \label{s-3}
    s_2^2-s_2+2cs_1-3c^2-2s_1s_3=0.
\end{equation}
Expressing $s_2$ from \eqref{s-1} we get
\begin{equation}\label{expr-s2}
 s_2=\frac{s_1^2-s_1+3c}{2}.   
\end{equation}
Next, using \eqref{s-1} and \eqref{s-2}, we can express $s_3$ in terms of $s_1$ and $c$ as
\begin{equation}\label{expr-s3}
   s_3=\frac{s_1^3-2s_1^2+7cs_1-s_1+3c}{6}.
\end{equation}
Substituting $s_2,s_3$ from \eqref{expr-s2} and \eqref{expr-s3} in \eqref{s-3} and factoring, we obtain
\begin{equation}
    \label{main-1}
    (s_1^2-3s_1+9c)(s_1^2+s_1+c+2)=0.
\end{equation}
This factorization appears naturally since we expect the triples $(p_1,p_1,p_1)$ and $(p_2,p_2,p_2)$ to be among the solutions of the system \eqref{eq1}-\eqref{eq3}, where $p_1,p_2$ are the fixed points of the map $f$; these fixed points exist for $c\leq \dfrac{1}{4}$, and are given by
$$p_{1}=\frac{1- \sqrt{1-4c}}{2}~~\text{and}~~ p_{2}=\frac{1+ \sqrt{1-4c}}{2}.$$
The two triples aforementioned correspond to the values $\dfrac{3(1\pm\sqrt{1-4c)}}{2}$ for $s_1$, which are indeed the solutions of the quadratic equation $s_1^2-3s_1+9c=0$ appearing from \eqref{main-1}. Furthermore, computing $s_2$ and $s_3$ for these values of $s_1$ using \eqref{expr-s2} and \eqref{expr-s3}, and then forming the cubic equations from the symmetric polynomials $s_1,s_2,s_3$, we obtain the following two equations:
$$(x-p_{1})^3=0 ~~~\text{and} ~~~(x-p_{2})^3=0.$$
It follows that the case $s_1^2-3s_1+9c=0$ that results from \eqref{main-1} leads to solutions of the system \eqref{eq1}-\eqref{eq3} which are not 3-cycles, since their components are equal. This means that the existence of 3-cycles may stem only from the second quadratic expression in \eqref{main-1}, that is
\begin{equation}
    \label{quadr}
    s_1^2+s_1+c+2=0.
\end{equation}
This equation admits real solutions if and only if its discriminant $\Delta=-4c-7$ is non-negative, that is,
\begin{equation}
    \label{cond-c}
    c\leq -\frac{7}{4}.
\end{equation}
To be completely certain of the existence of a 3-cycle, we must prove that the solutions $s_1$ of the quadratic equation \eqref{quadr} lead to corresponding values of $s_2$ and $s_3$ such that the cubic equations formed via Vi\`{e}te's formulas, with $x_1,x_2,x_3$ as unknowns, have three real (and distinct) roots.

The solutions of \eqref{quadr}, for $c\leq -\dfrac{7}{4}$, are
\begin{equation}
    \label{sol-s1}
    \widetilde{s_1}=\frac{-1- \sqrt{-4c-7}}{2}\,\,\,\,\,\text{and}\,\,\,\,\,\doubletilde{s_1}=\frac{-1+ \sqrt{-4c-7}}{2}.
\end{equation}
The computations below of the corresponding values of $s_2$ and $s_3$ can be simplified by successively reducing the degree of the variable $s_1$ in their expressions. We note that, given the equation \eqref{quadr}, and using \eqref{expr-s2} and \eqref{expr-s3},  we also have
\begin{equation*}
    s_2=-s_1+c-1\,\,\,\,\,\text{and}\,\,\,\,\,s_3=cs_1+c+1.
\end{equation*}
Thus, the values for $s_2$ obtained from \eqref{expr-s2} that correspond to the values of $s_1$ in \eqref{sol-s1} are:
\begin{equation}
    \label{sol-s2}
    \widetilde{s_2}=\frac{2c-1+\sqrt{-4c-7}}{2}\,\,\,\,\,\text{and}\,\,\,\,\,\doubletilde{s_2}=\frac{2c-1-\sqrt{-4c-7}}{2},
\end{equation}
while the values for $s_3$ obtained from \eqref{expr-s3} that correspond to the values of $s_1$ in \eqref{sol-s1} are:
\begin{equation}
    \label{sol-s3}
    \widetilde{s_3}=\frac{c+2-c\sqrt{-4c-7}}{2}\,\,\,\,\,\text{and}\,\,\,\,\,\doubletilde{s_3}=\frac{c+2+c\sqrt{-4c-7}}{2}.
\end{equation}
Using \eqref{sol-s1}, \eqref{sol-s2} and \eqref{sol-s3} and some straightforward algebraic simplifications, the discriminants\footnote{For a general cubic equation $Ax^3 + Bx^2 + Cx + D = 0$, the discriminant is given by the following expression: $B^2C^2 - 4AC^3 - 4B^3D - 27A^2D^2 + 18ABCD$.} of the two cubic equations that ensue are
\begin{equation}
    \label{discr-1}
    \Delta_1=16c^2-4c-7-8c\sqrt{-4c-7}\,\,\,\,\,\text{and}\,\,\,\,\,\Delta_2=16c^2-4c-7+8c\sqrt{-4c-7}.
\end{equation}

It is now easy to see that the discriminants in \eqref{discr-1} are positive since, for $c\leq -\dfrac{7}{4}$, we have $16c^2-4c-7>0$ and
$$(16c^2-4c-7)^2-(8c\sqrt{-4c-7})^2=(16c^2+4c+7)^2>0.$$

Therefore, for $c=-\dfrac{7}{4}$ we obtain one 3-cycle and for $c< -\dfrac{7}{4}$ we obtain two distinct 3-cycles for the map \eqref{map}.

In what follows, we discuss the stability of the 3-cycles using the criteria involving the first derivative (for hyperbolic cycles) and the second derivative or the Schwarzian derivative (for non-hyperbolic cycles); see, for example, Elaydi's book \cite[Sections 1.5 and 1.6]{El}.

Since for a 3-cycle $(x_1,x_2,x_3)$ of \eqref{map} we have
\begin{equation*}
    \label{prod-deriv}
    f'(x_1)f'(x_2)f'(x_3)=8x_1x_2x_3=8s_3,
\end{equation*}
we can easily deduce the stability of the 3-cycles using the first derivative test: if $|8s_3|<1$, then the cycle is asymptotically stable and if $|8s_3|>1$, then the cycle is unstable; the case $|8s_3|=1$ (non-hyperbolic cycle) is undecidable with the first derivative test. 

If $c<-\dfrac{7}{4}$, it is straightforward that
$$8\widetilde{s_3}-1=4c+7-4c\sqrt{-4c-7}=-\sqrt{-4c-7}\left(\sqrt{-4c-7}+4c\right)>0,$$
thus $8\widetilde{s_3}>1$ and the 3-cycle corresponding to $\widetilde{s_1}$ is unstable.\\
Similarly, we observe that
$$8\doubletilde{s_3}-1=4c+7+4c\sqrt{-4c-7}=\sqrt{-4c-7}\left(-\sqrt{-4c-7}+4c\right)<0,$$
therefore $8\doubletilde{s_3}<1$. The stability of the 3-cycle corresponding to $\doubletilde{s_1}$ is then guaranteed when $8\doubletilde{s_3}>-1$ or, equivalently,
\begin{equation*}
    \label{cond-stab}
    4c+9+4c\sqrt{-4c-7}>0,
\end{equation*}
which is in turn equivalent to
\begin{equation*}
   c\geq -\frac{9}{4}~~\hbox{and}~~~ (4c+9)^2-(4c\sqrt{-4c-7})^2>0
\end{equation*}
or
\begin{equation*}
    \label{cond-stab-2}
    c\geq -\frac{9}{4}~~\hbox{and}~~~64 c^{3}+128 c^{2}+72 c+81>0.
\end{equation*}
The discriminant of the above cubic polynomial is negative, therefore the polynomial has a unique real root, which is approximately $-1.768529$; in exact form, this root is 
\begin{equation*}
    \label{c-stab}
    \tilde{c}=-\frac{(7660 + 540\cdot \sqrt{201})^{1/3}}{24} - \frac{(7660 - 540\cdot \sqrt{201})^{1/3}}{24} - \frac{2}{3}.
\end{equation*}
It follows that the 3-cycle of \eqref{map} that corresponds to $\doubletilde{s_1}$ is asymptotically stable if $c\in \left(\tilde{c},-\frac{7}{4}\right)$.

For the special case $c=-\dfrac{7}{4}$, since we have $f'(x_1)f'(x_2)f'(x_3)=1$, we can try to use the stability test using the second derivative. Namely, we analyze the second derivative of the map $F:=f^{(3)}=f\circ f\circ f$ computed in the components of the 3-cycle $(x_1,x_2,x_3)$ which is formed, in this case, with the roots of the cubic equation
\begin{equation}
        \label{0}
        x^3+\frac{1}{2}x^2-\frac{9}{4}x-\frac{1}{8}=0;
    \end{equation}
see also \cite[p. 307]{HK}. One can easily check that $x_1\in (-2,-1),~x_2\in (1,2)$ and $x_3\in \left(-\frac{1}{10},0\right)$. Recall that, in this case, we also have $8x_1x_2x_3=1$.

We observe next that
\[
\begin{aligned}
F''(x)
&= f''(f(f(x)))\, [f'(f(x))]^2\, [f'(x)]^2 \\
&\quad + f'(f(f(x)))\, f''(f(x))\, [f'(x)]^2 \\
&\quad + f'(f(f(x)))\, f'(f(x))\, f''(x).
\end{aligned}
\]
Thus,
\begin{equation*}
F''(x_1)=32x_1^2x_2^2+16x_1^2x_3+8x_2x_3=32x_1^2x_2^2+2\cdot\frac{x_1}{x_2}+\frac{1}{x_1}=\frac{\dfrac{1}{16x_3^3}+2x_1^2+x_2}{x_1x_2},
\end{equation*}
and the analogous expressions for $F''(x_2)$ and $F''(x_3)$ obtained by circularly permuting the indices. Note now that 
$$\dfrac{1}{16x_3^3}+2x_1^2+x_2<-\frac{10^3}{16}+8+2<0.$$
Thus, $F''(x_1)>0$. Next, we compute
$$F''(x_2)=\frac{\dfrac{1}{16x_1^3}+2x_2^2+x_3}{x_2x_3},$$
and observe that
$$\dfrac{1}{16x_1^3}+2x_2^2+x_3>-\frac{1}{16}+2-\frac{1}{10}>0.$$
It follows from here that $F''(x_2)<0$. Finally, we have
$$F''(x_3)=\frac{\dfrac{1}{16x_2^3}+2x_3^2+x_1}{x_3x_1}$$
and
$$\dfrac{1}{16x_2^3}+2x_3^2+x_1<\frac{1}{16}+2\cdot\frac{1}{100}-1<0.$$
Therefore, $F''(x_3)<0$. Since $F''(x_1)$, $F''(x_2)$, $F''(x_3)$ are all non-zero, by the second derivative test for non-hyperbolic cycles, it follows that the 3-cycle in this case is unstable. The approximate values of this 3-cycle's components are $x_1\approx -1.746,~x_2\approx 1.301,~x_3\approx -0.054$.  

For the special case $c=\tilde c$, since we have $f'(x_1)f'(x_2)f'(x_3)=-1$, we can try to use the stability test using the Schwarzian derivative, namely by analyzing the sign of the Schwarzian of the map $F$ computed in the components of the 3-cycle $(x_1,x_2,x_3)$. Using the well-known formula
\begin{equation}\nonumber
\mathcal{S}\left(f_1\circ f_2\right)=\left(\left(\mathcal{S}f_1\right)\circ f_2\right)(f_2')^2+\mathcal{S}f_2,
\end{equation}
we have
\begin{equation*}
    \label{schwarz}
\mathcal{S}F(x_1)=\mathcal{S}f(x_3)\left(f'(x_1)f'(x_2)\right)^2+\mathcal{S}f(x_2)\left(f'(x_1)\right)^2+\mathcal{S}f(x_1),
\end{equation*}
and the other two analogous equalities.

Since $\mathcal{S}f(x)=-\displaystyle\frac{3}{2x^2}$ for $x\ne 0$, we obtain that $\mathcal{S}F(x_1), \mathcal{S}F(x_2),\mathcal{S}F(x_3)$ are negative, therefore, for $c=\tilde c$, the 3-cycle that corresponds to $\doubletilde{s_1}$ is still asymptotically stable.

\bigskip\smallskip

We can summarize the discussion above as follows:
\begin{enumerate}[{\bf (A)}]
    \item If $c>-\frac{7}{4}$, then there are no 3-cycles for the map \eqref{map};
    \item If $c=-\frac{7}{4}$, there is exactly one 3-cycle, having as components the roots of the cubic equation \eqref{0} and this cycle is unstable;
    \item If $c\in \left[\tilde c,-\frac{7}{4}\right)$, then there are two 3-cycles. One of the 3-cycles has as components the roots of the cubic equation
    \begin{equation}
        \label{1}
        x^3-\frac{-1-\sqrt{-4c-7}}{2}x^2+\frac{2c-1+\sqrt{-4c-7}}{2}x-\frac{c+2-c\sqrt{-4c-7}}{2}=0
    \end{equation}
    and is unstable. The other 3-cycle has as components the roots of the cubic equation
    \begin{equation}
        \label{2}
        x^3-\frac{-1+\sqrt{-4c-7}}{2}x^2+\frac{2c-1-\sqrt{-4c-7}}{2}x-\frac{c+2+c\sqrt{-4c-7}}{2}=0
    \end{equation}
    and is asymptotically stable;
    \item If $c<\tilde c$, then there are two 3-cycles for the map \eqref{map}, given by the roots of \eqref{1} and \eqref{2}, both 3-cycles being unstable.
\end{enumerate}

\begin{figure}[h]
\centering
\includegraphics[scale=0.65]{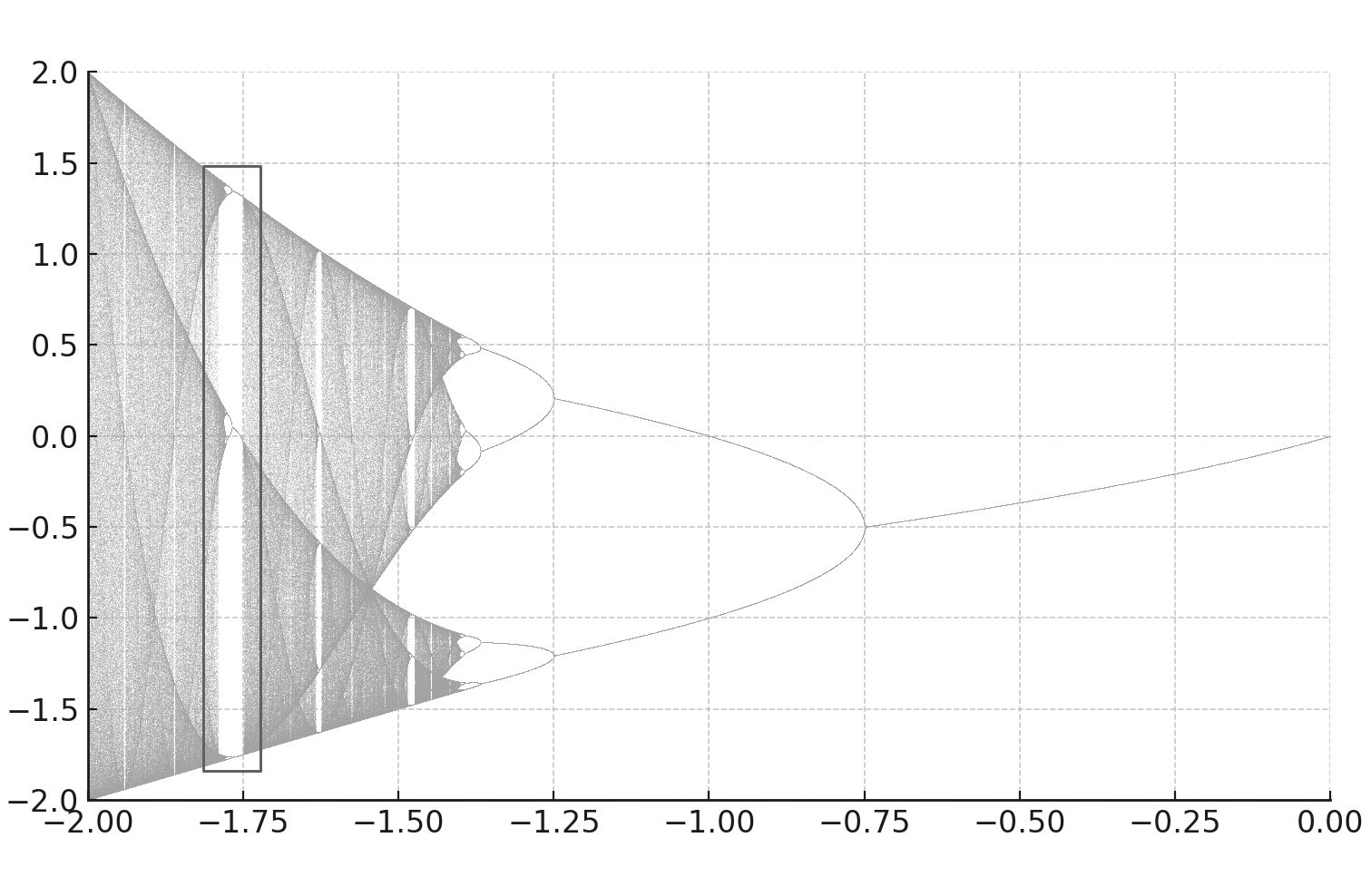}
\vspace*{-6mm}
\caption{Bifurcation diagram for the map \eqref{map} with parameter values $c\in [-2,0]$}
\label{fig-1}
\end{figure}

\begin{figure}[hb]
\centering
\includegraphics[scale=0.6]{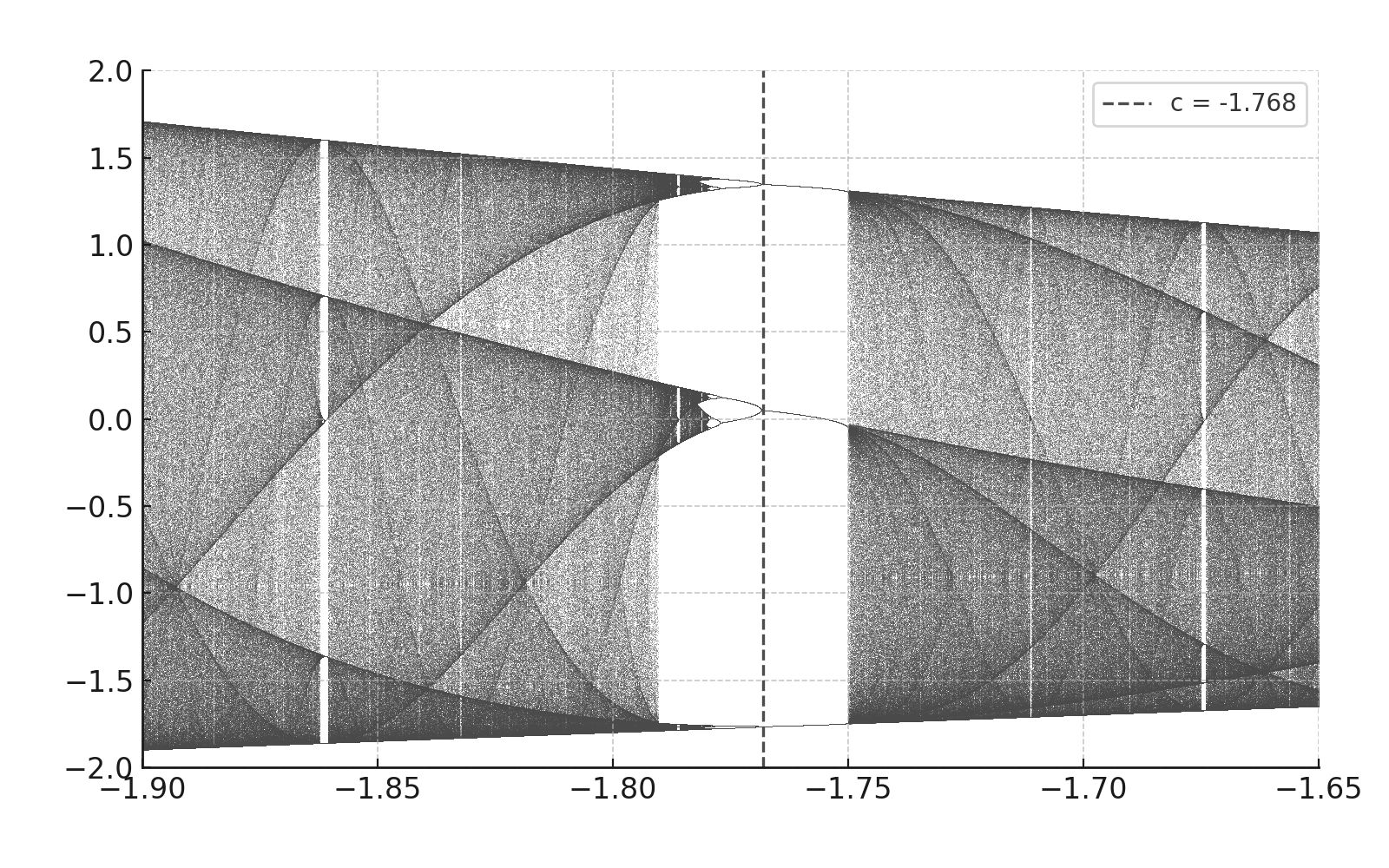}
\vspace*{-6mm}
\caption{Detail of the bifurcation diagram for the map \eqref{map}; the threshold value $\tilde c$, when stability of the 3-cycle vanishes, is marked with a dashed line}
\label{fig-2}
\end{figure}

We end this section with two figures that illustrate our analysis. In Figure \ref{fig-1}, the entire bifurcation diagram for the quadratic map \eqref{map} is presented for parameter values $c\in [-2,0]$. In Figure \ref{fig-2}, a detail of this diagram is shown; the 3-cycle appears at the right side of the delimited white ``stability oasis'', starting at $c=-\frac{7}{4}=-1.75$, but it rapidly fades, as the components of the stable 3-cycle begin a well-known period-doubling cascade (at $c=\tilde c$), which very rapidly descends back into chaos.
\vspace*{-0.4cm}

\section{The rebirth of 3-cycle for the logistic map}

It is well known that our quadratic map \eqref{map} is {\it conjugated} 
with the {\it logistic map}, given in terms of the parameter $r\in\mathbb R$ by 
\begin{equation}
\label{logistic}
g_r:\mathbb{R}\to \mathbb{R},~~g_r(y)=ry(1-y);
\end{equation}
see, for example, \cite[p. 61]{HK} or \cite[p. 48]{El}. Usually, taking into account the origin of the logistic map and its parameter in modeling biological systems, one restricts the values of the logistic parameter $r$ to belong to the interval $(0,4]$. As before, we will suppress the parameter $r$ as a subscript and simply denote the logistic map by $g$.

The works of Saha--Strogatz \cite{SS}, Bechhoefer \cite{Be}, and Gordon \cite{Go} offer elementary arguments for determining the smallest value of the logistic parameter \( r \) at which a 3-cycle appears and chaos consequently emerges in the logistic map. Despite their apparent simplicity, each of these arguments involves a subtle ``twist'': a clever change of variables in \cite{SS} and \cite{Be}, and an implicit use of Carvalho’s lemma \cite[p.~41]{El} in \cite{Go}. More recent proofs were given by Burm-Fishback \cite{BF}, Calvis \cite{Ca} and  Zhang \cite{Zh}.

For a given parameter $r$, the logistic map \eqref{logistic} is conjugated with the quadratic map \eqref{map} with parameter
\begin{equation}
\label{r-c}   
c=-\frac{r(r-2)}{4}.
\end{equation}
More precisely, the conjugation map is given by
\begin{equation*}
\label{conj-map}
h:\mathbb{R}\to \mathbb{R}, ~~h(x)=-rx+\frac{r}{2}
\end{equation*}
and it verifies the conjugacy relation
\begin{equation}
\label{conj-def}
f\circ h=h\circ g.
\end{equation}
Using the definition of a 3-cycle and property \eqref{conj-def}, it is easy to prove that
\begin{equation*}
    \label{equiv}
    (y_1,y_2,y_3)\hbox{ is a 3-cycle for }g \Longleftrightarrow (h(y_1),h(y_2),h(y_3))\hbox{ is a 3-cycle for }f.
\end{equation*}
It follows that a 3-cycle exists for the logistic map $g$ if and only if the corresponding value of the parameter $c$ for the quadratic map \eqref{map} stated in \eqref{r-c} verifies the condition \eqref{cond-c}, namely
$$-\frac{r(r-2)}{4}\leq -\frac{7}{4},$$
or, more precisely,
\begin{equation*}
    \label{cond-r}
    r \in \left(-\infty, 1-2\sqrt{2}\right] \cup \left[1+2\sqrt{2},\infty\right).
\end{equation*}
The condition \( r \geq 1 + 2\sqrt{2} \approx 3.828427 \) is widely recognized in the literature and has been proved in various ways. Note also that $1-2\sqrt{2}\approx -1.828427$.

Moreover, for a 3-cycle $(y_1,y_2y_3)$ of the logistic map \eqref{logistic} we have
\begin{equation*}
    \label{connection-stab}
    g'(y_1)g'(y_2)g'(y_3)=f'(h(y_1))f'(h(y_2))f'(h(y_3)),
\end{equation*}
meaning that the stability (studied with the first derivative test) of the 3-cycle $(y_1,y_2y_3)$ is equivalent with the stability of the 3-cycle $(h(y_1),h(y_2),h(y_3))$ of \eqref{map}. This leads to the condition
\begin{equation*}
    \label{r-max}
    c=-\frac{r(r-2)}{4}> \tilde c,
\end{equation*}
or, equivalently,
$$r\in \left(r_{min},r_{max}\right)$$
with
$$r_{min}=1-\frac{\sqrt{132+6 \left(7660+540 \sqrt{201}\right)^{{1}/{3}}+6 \left(7660-540 \sqrt{201}\right)^{{1}/{3}}}}{6}\approx -1.841499$$
and
$$r_{max}=1+\frac{\sqrt{132+6 \left(7660+540 \sqrt{201}\right)^{{1}/{3}}+6 \left(7660-540 \sqrt{201}\right)^{{1}/{3}}}}{6}\approx 3.841499.$$
The value $r_{max}$ can be re-written as
$$r_{max}=1+\sqrt{\frac{11}{3}+\left(\frac{1915}{54}+\frac{5 \sqrt{201}}{2}\right)^{1 / 3}+\left(\frac{1915}{54}-\frac{5 \sqrt{201}}{2}\right)^{1 / 3}},$$
which is the maximum value of the logistic parameter $r$ that supports a stable 3-cycle given in \cite{Go}.

Therefore, for the logistic map, without extra assumptions on the parameter $r$, a 3-cycle exists and is asymptotically stable if $r\in \left(r_{min},1-2\sqrt{2}\right)$ or $r\in \left(1+2\sqrt{2},r_{max}\right)$.

\subsection*{Acknowledgments}
The first author was supported by an AMS-Simons Research Enhancement Grant for PUI Faculty.

\end{document}